\documentclass{amsart}

\usepackage{amsmath}
\usepackage{amssymb}

\newtheorem{theorem}{Theorem}[section]

\newtheorem{lemma}[theorem]{Lemma}
\newtheorem{proposition}[theorem]{Proposition}

\begin{document}

\title{The Kaplansky Test Problems for $\aleph _1$-Separable Groups}

\author{Paul C. Eklof}
\address{Department of Mathematics,University of California, Irvine 92697}
\email{pceklof@uci.edu}
\thanks{Travel supported by NSF Grant DMS-9501415.} 

\author{Saharon Shelah}
\address{Department of Mathematics, University of Wisconsin, WI}
\curraddr{Institute of Mathematics, Hebrew University, Jerusalem 91904}
\email{shelah@math.huji.ac.il}
\thanks{Research supported by German-Israeli Foundation for Scientific Research \&
Development Grant No. G-294.081.06/93. Pub. No. 625}


\date{}


\begin{abstract}
We answer a long-standing open question by proving in ordinary set theory,
ZFC, that the Kaplansky test problems have negative answers for $\aleph _1$%
-separable abelian groups of cardinality $\aleph _1$. In fact, there is an $%
\aleph _1$-separable abelian group $M$ such that $M$ is isomorphic to $%
M\oplus M\oplus M$ but not to $M\oplus M$. We also derive some relevant
information about the endomorphism ring of $M$.
\end{abstract}
\maketitle

\section*{Introduction}

Kaplansky \cite[pp. 12f]{K} posed two test problems in order to ``know when
we have a satisfactory [structure] theorem. ... We suggest that a tangible
criterion be employed: the success of the alleged structure theorem in
solving an explicit problem.'' The two problems were:

\begin{quote}
(I) If $A$ is isomorphic to a direct summand of $B$ and conversely, are $A$
and $B$ isomorphic?

(II) If $A\oplus A$ and $B\oplus B$ are isomorphic, are $A$ and $B$
isomorphic?
\end{quote}

In fact, he says (\cite[p. 75]{K}) that he invented the problems ``to show
that Ulm's theorem [a structure theory for countable abelian $p$-groups]
could really be used''. For some other classes of abelian groups, such as
finitely-generated groups, free groups, divisible groups, or completely
decomposable torsion-free groups, the existence of a structure theory leads
to an affirmative answer to the test problems. On the other hand, negative
answers are taken as evidence of the absence of a useful classification
theorem for a given class; Kaplansky says ``I believe their defeat is
convincing evidence that no reasonable invariants exist'' \cite[p. 75]{K}.
Negative answers to both questions have been proven, for example, for the
class of uncountable abelian $p$-groups and for the class of countable
torsion-free abelian groups.

Of particular interest is the method developed by Corner (cf. \cite{C1}, 
\cite{C2},\cite{CG}) which, by realizing certain rings as endomorphism rings
of groups, provides negative answers to both test problems (for a given
class) as special cases of an even more extreme pathology. More precisely,
Corner's method --- where applicable --- yields, for any positive integer $r 
$, an abelian group $G_r$ (in the class) such that for any positive integers 
$m$ and $k$, the direct sum of $m$ copies of $G_r$ is isomorphic to the
direct sum of $k$ copies of $G_r$ if and only if $m$ is congruent to $k$ mod 
$r$. (See, for example, \cite{C2} or \cite[Thm 91.6, p. 145]{F2}.) Then we
obtain negative answers to both test problems by letting $A=G_2$ ($\cong
G_2\oplus G_2\oplus G_2$) and $B=G_2\oplus G_2$.

Our focus here is on the class of $\aleph _1$-separable abelian groups (of
cardinality $\aleph _1$). We will prove, in ordinary set theory (ZFC), that
both test problems have negative answers by deriving the Corner pathology:

\begin{theorem}
\label{main}For any positive integer $r$ there is an $\aleph _1$-separable
group $M=M_r$ of cardinality $\aleph _1$ such that for any positive integers 
$m$ and $k$, $M^m$ is isomorphic to $M^k$ if and only if $m$ is congruent to 
$k$ mod $r$.
\end{theorem}

\noindent (Here $M^m$ denotes the direct sum of $m$ copies of $M$.) We do
not determine the endomorphism ring of $M$, even modulo an ideal. However,
we can derive a property of the endomorphism ring of $M$ which is sufficient
to imply the Corner pathology: see section \ref{endo}.

\smallskip\ 

A group $M$ is called $\aleph _{1}$-separable \cite[p. 184]{F}
(respectively, strongly $\aleph _{1}$-free) if it is abelian and every
countable subset is contained in a countable free direct summand of $M$
(resp., contained in a countable free subgroup $H$ which is a direct summand
of every countable subgroup of $M$ containing $H$). Obviously, an $\aleph
_{1}$-separable group is strongly $\aleph _{1}$-free, so a negative answer
to one of the test problems for the class of $\aleph _{1}$-separable groups
implies a negative answer to the problem for the class of strongly $\aleph
_{1}$-free groups. (It is independent of ZFC whether these classes are
different for groups of cardinality $\aleph _{1}$: the weak Continuum
Hypothesis ($2^{\aleph _{0}}<2^{\aleph _{1}}$) implies that there are
strongly $\aleph _{1}$-free groups of cardinality $\aleph _{1}$ which are
not $\aleph _{1}$-separable; on the other hand, Martin's Axiom (MA) plus the
negation of the Continuum Hypothesis ($\lnot $CH) implies that every
strongly $\aleph _{1}$-free group of cardinality $\aleph _{1}$ is $\aleph
_{1}$-separable; cf.\cite{M80} )

Dugas and G\"{o}bel \cite{DGa} proved that ZFC + $2^{\aleph _{0}}<2^{\aleph
_{1}}$ implies that the Corner pathology exists for the class of strongly $%
\aleph _{1}$-free groups of cardinality $\aleph _{1}$; in fact, they showed
that there is a strongly $\aleph _{1}$-free group $G$ whose endomorphism
ring is an appropriate ring (the ring $A=$ $A_{r}$ of the next section).
(See also \cite{G}.) This group $G$ cannot be $\aleph _{1}$-separable since
the endomorphism ring of an $\aleph _{1}$-separable group has too many
idempotents. However, Thom\'{e} (\cite{Th1} and \cite{Th2}) showed that ZFC
plus V = L (G\"{o}del's Axiom of Constructibility) implies the Corner
pathology for $\aleph _{1}$-separable groups of cardinality $\aleph _{1}$;
he did this by constructing an $\aleph _{1}$-separable $G$ such that $%
\operatorname{End}(G)$ is a split extension of $A$ by $I$ (in the sense of 
\cite[p. 277]{C3}), where $I$ is the ideal of endomorphisms with a countable
image.

It follows from known structure theorems for the class of $\aleph _1$%
-separable groups of cardinality $\aleph _1$ under the hypothesis MA + $\neg 
$CH that the Dugas-G\"obel and Thom\'e realization results are \textit{not}
theorems of ZFC (cf. \cite{E83a} or \cite{M83}). The fact that there \textit{%
are} positive structure theorems for the class of $\aleph _1$-separable
groups assuming MA + $\neg $CH or the stronger Proper Forcing Axiom (PFA)
--- see, for example, \cite{E83b} or \cite{M87} --- led to the question of
whether the Kaplansky test problems could have affirmative answers for this
class assuming, say, PFA. Thom\'e \cite{Th2} gave a negative answer to the
second test problem in ZFC, using a result of J\'onsson \cite{Jon} for
countable torsion-free groups; however, till now, the first test problem as
well as the Corner pathology were open (in ZFC).

\smallskip\ 

Our construction of the Corner pathology involves a direct construction of
the pathological group $M$ using a tree-like ladder system and a ``countable
template'' which comes from the Corner example for countable torsion-free
groups. A key role is played by a paper of G\"{o}bel and Goldsmith \cite{GG}
which --- while it does not itself prove any new results about the Kaplansky
test problems for strongly $\aleph _{1}$-free or $\aleph _{1}$-separable
groups --- provides the tools for creating a suitable template from the
Corner example.

\section{The countable template}

Fix a positive integer $r$. For this $r$, let $A=A_{r}$ be the countable
ring constructed by Corner in \cite{C2}. (See also \cite[p. 146]{F2}.)
Specifically, $A$ is the ring freely generated by symbols $\rho _{i}$ and $%
\sigma _{i}$ ($i=0,1,...,r$) subject to the relations 
\begin{equation*}
\rho _{j}\sigma _{i}=\left\{ 
\begin{array}{ll}
1 & \text{if }i=j \\ 
0 & \text{otherwise}
\end{array}
\right. 
\end{equation*}
\noindent and
\begin{equation*}
\sum_{i=0}^{r}\sigma _{i}\rho _{i}=1\text{.}
\end{equation*}
Then $A$ is free as an abelian group, and $\sigma _{0}\rho _{0},...,\sigma
_{r}\rho _{r}$ are pairwise orthogonal idempotents. Moreover, if $M$ is a
right $A$-module, then $M=M\sigma _{0}\rho _{0}\oplus M\sigma _{1}\rho
_{1}\oplus ...\oplus M\sigma _{r}\rho _{r}$ and $M\sigma _{i}\rho _{i}\cong M
$ because $\sigma _{i}\rho _{i}\sigma _{i}:M\rightarrow M\sigma _{i}\rho _{i}
$ and $\rho _{i}\sigma _{i}\rho _{i}:M\sigma _{i}\rho _{i}\rightarrow M$ are
inverses; therefore $M\cong M^{r+1}$. 

Our construction will work for any
countable torsion-free ring $A$ whose additive subgroup is free; but
hereafter $A$ will denote the ring $A_{r}$ just defined.

Corner shows that there is a torsion-free countable abelian group $G$ whose
endomorphism ring is $A$; thus $G$ is an  $A$-module and hence $G^{{}}\cong
G^{r+1}$. Furthermore, he shows that $G^{\ell }$ is not isomorphic to $G^{n}$
if $1\leq \ell <n\leq r$, and hence $G^{m}$ is not isomorphic to $G^{k}$ if $%
m$ is not congruent to $k$ mod $r$. We shall require these and further
properties of $G$, which we summarize in the following:

\begin{proposition}
\label{tem}There are countable free $A$-modules $B\subseteq H$ such that $%
G\cong H/B$ and $B$ is the union of a chain of free $A$-modules, $%
B=\bigcup_{n\in \omega }B_n$, such that $B_0=0$ and for all $n\in \omega $, $%
H/B_n$ and $B_{n+1}/B_n$ are free $A$-modules of rank $\omega $. Moreover
for any positive integers $m$ and $k$, if $m$ is not congruent to $k$ mod $r$%
, then $G^m\oplus {\bf Z}^{(\omega )}$ is not isomorphic to $G^k\oplus {\bf Z%
}^{(\omega )}$.
\end{proposition}

The main work in proving Proposition \ref{tem} will be done in two lemmas
from \cite{GG}. For the first one, we give a revised proof (cf. \cite[p. 343]
{GG}). We maintain the notation above.

\begin{lemma}
\label{gg}The group $G$ is the union, $G=\bigcup_{n\geq 1}G_n$, of an
increasing chain of free $A$-modules.
\end{lemma}

\begin{proof} By \cite[p. 699]{C1} $G$ is the pure closure $\left\langle
G_1\right\rangle _{*}$ in $\hat A$ of a free $A$-module $G_1=\bigoplus_{i\in
I}e_iA\oplus A$ containing $A$. Here $\hat A$ is the natural, or $\mathbf{Z}$%
-adic, completion of $A$ (cf. \cite[p. 692]{C1}). We will define inductively 
$G_n=\bigoplus_{i\in I}e_{i,n}A\oplus A$ such that $G_n\supseteq G_{n-1}$
and for all $i\in I$, $ne_{i,n}+A=e_{i,n-1}+A$. Let $e_{i,1}=e_i$ for all $%
i\in I$. If $G_{n-1}\subseteq G$ has been defined for some $n>1$, then since 
$A$ is dense in $\hat A$, there exists $e_{i,n}\in \hat A$ such that $%
ne_{i,n}+A=e_{i,n-1}+A$; say $ne_{i,n}=e_{i,n-1}+a_i$. By the definition of $%
G$, $e_{i,n}\in G$. We need to show that $\{e_{i,n}:i\in I\}\cup \{1\}$ is $%
A $-linearly independent. Suppose that $\Sigma _{i\in I}e_{i,n}c_i+1\cdot
c_0=0 $ for some $c_0,c_i\in A$. Then $\Sigma _{i\in I}ne_{i,n}c_i+nc_0=0$,
so $\Sigma _{i\in I}e_{i,n-1}c_i+1\cdot (\Sigma _{i\in I}a_ic_i+nc_0)=0$. By
the $A$-linear independence of $\{e_{i,n-1}:i\in I\}\cup \{1\}$, we can
conclude that each $c_i$ equals $0$ and hence also $c_0$ equals $0$. This
completes the definition of $G_n$.

It remains to prove that $G\subseteq \bigcup_{n\geq 1}G_n$. Let $g\in
G\setminus G_1$. For some $n>1$, $ng\in G_1$. We claim that $g\in G_n$.
Since $ng\in G_{n-1}$, $ng=\Sigma _{i\in I}e_{i,n-1}c_i+c_0$ for some $%
c_i,c_0\in A$. Then 
\begin{equation*}
ng=\Sigma _{i\in I}(ne_{i,n}-a_i)c_i+c_0=n\Sigma _{i\in
I}e_{i,n}c_i+a^{\prime } 
\end{equation*}
for some $a^{\prime }\in A$. Since $A$ is pure in $\hat A$, $a^{\prime
}=na^{\prime \prime }$ for some $a^{\prime \prime }\in A$. Thus $g=\Sigma
_{i\in I}e_{i,n}c_i+a^{\prime \prime }\in G_n$.   
\end{proof}

\smallskip

The second lemma is proved in \cite[Lemma 2.5]{GG} generalizing a result in 
\cite[Lemma XII.1.4]{EM}. We state it here for the sake of completeness.

\begin{lemma}
\label{gg2} Let  $G$ be a
countable $A$-module which is the union, $G=\bigcup_{n\geq 1}G_n$, of an
increasing chain of free $A$-modules, then there exist countable free $A$%
-modules $B\subseteq H$ such that $G\cong H/B$ and $B$ is the union of a
chain of free $A$-modules, $B=\bigcup_{n\geq 1}B_n$, such that for all   $%
n\geq 1$, $H/B_n$ and $B_{n+1}/B_n$ are free $A$-modules.  $\square$ 

\end{lemma}

\textsc{proof of Proposition \ref{tem}. }The existence of $H$, $B$, and the $%
B_{n}$ is now an immediate consequence of Lemmas \ref{gg} and \ref{gg2}. All
that is left to show is that if $m$ is not congruent to $k$ mod $r$, then $%
G^{m}\oplus \mathbf{Z}^{(\omega )}$ is not isomorphic to $G^{k}\oplus 
\mathbf{Z}^{(\omega )}$. Since $G^{m}$ is not isomorphic to $G^{k}$, it is
enough to show that $R_{\mathbf{Z}}(G^{l}\oplus \mathbf{Z}^{(\omega )})=G^{l}
$ for any $l\in \omega $. Here $R_{\mathbf{Z}}(N)$ is the $\mathbf{Z}$%
-radical of $N$, that is, $R_{\mathbf{Z}}(N)=\cap \{\ker (\varphi ):\varphi
:N\rightarrow \mathbf{Z}\}$. (See, for example, \cite[pp. 289f]{EM}.) To
show that $R_{\mathbf{Z}}(G^{l}\oplus \mathbf{Z}^{(\omega )})=G^{l}$ it is
enough to show that $\operatorname{Hom}(G^{l},\mathbf{Z})=0$, or, equivalently, $%
\operatorname{Hom}(G,\mathbf{Z})=0$. This follows from Observation 2.7 of \cite
{GG}, but we give here a self-contained argument based on the notation of
Lemma \ref{gg}. Suppose $\psi \in \operatorname{Hom}(G,\mathbf{Z})$; we can
regard $\psi $ as an endomorphism of $G$ by identifying $\mathbf{Z}$ with
the subgroup $\left\langle 1\right\rangle $ of $A\subseteq G$ which is
generated by the unit $1$ of $A$. Since the endomorphism ring of $G$ is $A$,
there is $a\in A$ such that $\psi (g)=ga$ for all $g\in G$. By considering $%
\psi (1)=1a=a$, we see that $a\in \left\langle 1\right\rangle $. Now
consider $\psi (e_{i})$ for any $e_{i}$; since $\psi (e_{i})=e_{i}a$ and
since $e_{i}A\cap \left\langle 1\right\rangle =\{0\}$ we see that $a=0$.  $\square$

\section{The main construction}

Fix a positive integer $r$ and let $A,H,B,B_n$ and $G$ be as in Proposition 
\ref{tem}. For each $n\in \omega $, fix a basis $\{b_{n,i}+B_n:i\in \omega
\} $ of $B_{n+1}/B_n$ (as $A$-module). Also, fix a set of representatives $%
\{h_i:i\in \omega \}$ for $H/B$ where $h_0=0$; thus each coset $h+B$ equals $%
h_i+B$ for a unique $i\in \omega $.

Fix a stationary subset $E$ of $\omega _1$ consisting of limit ordinals and
a ladder system $\{\eta _\delta :\delta \in E\}$. That is, for every $\delta 
$ in $E$, $\eta _\delta :\omega \rightarrow \delta $ is a strictly
increasing function whose range is cofinal in $\delta $; we shall also
choose $\eta _\delta $ so that its range is disjoint from $E$. Furthermore,
we choose a ladder system which is \textit{tree-like}, that is, for all $%
\delta ,\gamma \in E$ and $n,m\in \omega $, $\eta _\delta (n)=\eta _\gamma
(m)$ implies that $m=n$ and $\eta _\delta (l)=\eta _\gamma (l)$ for all $l<n$
(cf. \cite[pp. 368, 386]{EM}).

Inductively define free $A$-modules $M_\beta $ ($\beta <\omega _1$) as
follows: if $\beta $ is a limit ordinal, $M_\beta =\bigcup_{\alpha <\beta
}M_\alpha $; if $\beta =\alpha +1$ where $\alpha \notin E$, let
\begin{equation*}
M_\beta =M_\alpha \oplus \bigoplus_{i\in \omega }x_{\alpha ,i}A\text{.} 
\end{equation*}
If $\beta =\delta +1$ where $\delta \in E$, define an embedding $\iota
_\delta :B\rightarrow M_\delta $ by sending the basis element $b_{n,i}$ to $%
x_{\eta _\delta (n),i}$. Essentially $M_{\delta +1}$ will be defined to be
the pushout of 
\begin{equation*}
\begin{array}{lll}
M_\delta &  &  \\ 
\uparrow \iota _\delta &  &  \\ 
B & \hookrightarrow & H
\end{array}
\end{equation*}
but we will be more explicit in order to avoid the necessity of identifying
isomorphic copies. Let $y_{\delta ,0}=0$ and let $\{y_{\delta ,i}:i\in
\omega \setminus \{0\}\}$ be a new set of distinct elements (not in $%
M_\delta $). Then define $M_{\delta +1}$ to be $\{y_{\delta ,i}+u:u\in
M_\delta $, $i\in \omega \}$ where the operations on $M_{\delta +1}$ extend
those on $M_\delta $ and are otherwise determined by the rules 
\begin{equation*}
\begin{array}{lll}
y_{\delta ,i}+y_{\delta ,j}=y_{\delta ,k}+\iota _\delta (b) & \text{if} & 
h_i+h_j=h_k+b \\ 
y_{\delta ,i}a=y_{\delta ,\ell }+\iota _\delta (b) & \text{if} & h_ia=h_\ell
+b
\end{array}
\end{equation*}

\noindent where $b\in B$ and $a\in A$. Then there is an embedding $\theta
_\delta :H\rightarrow M_{\delta +1}$ extending $\iota _\delta $ which takes $%
h_i$ to $y_{\delta ,i}$ and induces an isomorphism of $H/B$ with $M_{\delta
+1}/M_\delta $.

This completes the inductive definition of the $M_\beta $. Let $%
M=\bigcup_{\beta <\omega _1}M_\beta $. Note that it follows from the
construction that every element of $M$ has a unique representation in the
form
\begin{equation*}
\sum_{j=1}^sy_{\delta _j,n_j}+\sum_{\ell =1}^tx_{\alpha _\ell ,i_\ell
}a_\ell 
\end{equation*}
where $\delta _1<\delta _2<...<\delta _s$ are elements of $E$, $n_j\in
\omega \setminus \{0\}$, $\alpha _\ell \in \omega _1\setminus E$, $i_\ell
\in \omega $, $a_\ell \in A$, and the pairs $(\alpha _\ell ,i_\ell )$ ($\ell
=1,...,t$) are distinct.

Since $M$ is constructed to be an $A$-module, $M$ is isomorphic to $M^{r+1}$%
. We claim that

\begin{quote}
($\dagger $) $M$ is $\aleph _1$-separable; in fact for all $\alpha <\omega
_1 $, $M_{\alpha +1}$ is a free direct summand of $M$.
\end{quote}

\noindent Assuming this for the moment, we can show that

\begin{quote}
($\dag \dag $) $M^m$ is not isomorphic to $M^k$ if $m$ is not congruent

to $k$ mod $r$.
\end{quote}

In brief this is because $M^m$ and $M^k$ are not quotient-equivalent (cf. 
\cite[pp. 251f]{EM}) since for all $\delta \in E$, $(M_{\delta +1}/M_\delta
)^m\oplus \mathbf{Z}^{(\omega )}$ is not isomorphic to $(M_{\delta
+1}/M_\delta )^k\oplus \mathbf{Z}^{(\omega )}$ by Proposition \ref{tem}. In
more detail, if there is an isomorphism $\varphi :M^m\rightarrow M^k$, then
there is a closed unbounded subset $C$ of $\omega _1$ such that for $\beta
\in C$, $\varphi [M_\beta ^m]=M_\beta ^k.$ Since $E$ is stationary in $%
\omega _1$, there exist $\delta \in C\cap E$; choose $\beta >\delta $ such
that $\beta \in C$. Then $\varphi $ induces an isomorphism of $M_\beta
^m/M_\delta ^m$ with $M_\beta ^k/M_\delta ^k$. Since $M_\beta /M_{\delta +1}$
is free (of infinite rank) by ($\dagger $), we can conclude that 
\begin{equation*}
\begin{array}{c}
(M_{\delta +1}/M_\delta )^m\oplus \mathbf{Z}^{(\omega )}\cong (M_{\delta
+1}^m/M_\delta ^m)\oplus (M_\beta ^m/M_{\delta +1}^m)\cong M_\beta
^m/M_\delta ^m\cong M_\beta ^k/M_\delta ^k \\ 
\cong (M_{\delta +1}^k/M_\delta ^k)\oplus (M_\beta ^k/M_{\delta +1}^k)\cong
(M_{\delta +1}/M_\delta )^k\oplus \mathbf{Z}^{(\omega )}
\end{array}
\end{equation*}
which contradicts Proposition \ref{tem}.

We are left with the task of proving ($\dagger $). First we shall show that
each $M_{\alpha +1}$ is a direct summand of $M$ by defining a projection $%
\pi _{\alpha }$ of $M$ onto $M_{\alpha +1}$ (that is, $\pi _{\alpha
}|M_{\alpha +1}$ is the identity). For every integer $k$ there is a
projection $\rho _{k}:H\rightarrow B_{k+1}$ since $H/B_{k+1}$ is free. Given 
$\alpha $, for each $\delta \in E$ with $\delta >\alpha $, let $k_{\delta }$
be the maximal integer $k$ such that $\eta _{\delta }(k)\leq \alpha $. For
each $\delta \in E$, we let $\pi _{\alpha }$ act like $\rho _{k_{\delta }}$
on the isomorphic copy, $\theta _{\delta }[H]$, of $H$. More precisely, for
each element $z$ of  $\theta _{\delta }[H]$, define $\pi _{\alpha }(z)$ to
be $\theta _{\delta }(\rho _{k_{\delta }}(\theta _{\delta }^{-1}(z)))$; if $%
\nu \notin $ $\bigcup \{\operatorname{ran}(\eta _{\delta }):\delta \in E\}$ and $%
\nu >\alpha $, define $\pi _{\alpha }(x_{\nu ,i})=0$. Extend to an arbitrary
element of $M$ by additivity; this will define a homomorphism on $M$
provided that $\pi _{\alpha }$ is well-defined. It is easy to see, using the
unique representation of elements, that the question of well-definition
reduces to showing that the definition of $\pi _{\alpha }(x_{\beta ,i})$ for 
$x_{\beta ,i}\in $ $\theta _{\delta }[H]$ is independent of $\delta $. If $%
\beta \leq \alpha $, then $\pi _{\alpha }(x_{\beta ,i})=x_{\beta ,i}$. Say $%
\beta >\alpha $ and $\beta =\eta _{\delta }(n)=\eta _{\gamma }(n)$; by the
tree-like property, $\eta _{\delta }(m)=\eta _{\gamma }(m)$ for all $m\leq n$%
, and hence $k_{\delta }=k_{\gamma }$. Hence $\pi _{\alpha }(x_{\beta ,i})$
is well-defined because $\rho _{k_{\delta }}=\rho _{k_{\gamma }}$ and thus $%
\theta _{\delta }(\rho _{k_{\delta }}(\theta _{\delta }^{-1}(x_{\beta
,i})))=\theta _{\gamma }(\rho _{k_{\gamma }}(\theta _{\gamma }^{-1}(x_{\beta
,i})))$.

It remains to prove that each $M_\beta $ is $\aleph _1$-free (as abelian
group). Since $A$ is free as abelian group, it suffices to show that $%
M_{\delta +1}$ is a free $A$-module for every $\delta \in E$. We will
inductively define $S_n$ so that 
\begin{equation*}
B=\bigcup_{n\in \omega }S_n\cup \{x_{\nu ,i}:\nu \in \delta \setminus (E\cup
\bigcup \{\operatorname{ran}(\eta _\mu ):\mu \in E\cap (\delta +1)\}),i\in \omega
\} 
\end{equation*}
is an $A$-basis of $M_{\delta +1}$. Let $S_0$ be the image under $\theta
_\delta $ of a basis of $H$. Fix a bijection $\psi :\omega \rightarrow E\cap
\delta $; also, for convenience, let $\psi (-1)=\delta $. Suppose that $S_m$
has been defined for $m\leq n$ so that $\bigcup_{m\leq n}S_m$ is $A$%
-linearly independent and generates $\bigcup \{\theta _{\psi (m)}[H]:-1\leq
m<n\}$. Let $\gamma =\psi (n)$ and let $k=k_n$ be maximal such that $\eta
_\gamma (k)=\eta _{\psi (m)}(k)$ for some $-1\leq m<n$. Notice that $%
\{x_{\eta _\gamma (\ell ),i}:\ell \leq k,i\in \omega \}$ is contained in the 
$A$-submodule generated by $\bigcup_{m\leq n}S_m$. Since $H/B_{k+1}$ is $A$%
-free, we can write $H=B_{k+1}\oplus C_k$ for some $A$-free module $C_k$ ($%
=\ker (\rho _k)$); let $S_{n+1}$ be the image under $\theta _\gamma $ of a
basis of $C_k$. This completes the inductive construction. One can then
easily verify that $B$ is an $A$-basis of $M_{\delta +1}$; indeed, the fact
that $\bigcup_{m\leq n}S_m$ is $A$-linearly independent can be proved by
induction on $n$, using the unique representation of elements of $M$ to show
that if $\sum_{i=1}^rz_ia_i\in \left\langle \bigcup_{m\leq
n}S_m\right\rangle $, where $z_1,...,z_r$ are distinct elements of $S_{n+1}$%
, then $a_i=0$ for all $i=1,...,r$.

\bigskip 

\section{The endomorphism ring of M\label{endo}}

While we cannot show that $\operatorname{End}(M)$ is a split extension of $A$ by
an ideal, we can obtain enough information about $\operatorname{End}(M)$ to imply
the negative results on the Kaplansky test problems. (A similar idea is used
in \cite[p. 118]{Sh381}.)

The ring $A$ is naturally a subring of $\operatorname{End}(M)$. We say that $A$
is \textit{algebraically closed} in $\operatorname{End}(M)$ when every finite set
of ring equations with parameters from $A$ (i.e., polynomials in several
variables over $A$) which is satisfied in $\operatorname{End}(M)$ is also
satisfied in $A$.

\begin{proposition}
\label{algclsd}
If $A=A_r$ is as in section 1, and $A$ is algebraically closed in $\operatorname{End}(M)$, then for any positive integers $m$ and $k$, $M^m$ is isomorphic to $M^k$ if and only if $m$ is congruent to $k$ mod $r$.
\end{proposition}

\begin{proof} Since $M$ is an $A$-module, $M\cong M^{r+1}.$ If $M^\ell \,$
is isomorphic to $M^n$ where $1\leq \ell <n\leq r$, then $\sum_{i=1}^\ell
M\sigma _i\rho _i\cong $ $\sum_{i=1}^nM\sigma _i\rho _i$, so by Lemma 2 of 
\cite{C2}, there are elements $x$ and $y$ of $\operatorname{End}(M)$ such that $%
xy=\sum_{i=1}^\ell \sigma _i\rho _i$ and $yx=$ $\sum_{i=1}^n\sigma _i\rho _i$%
. So by hypothesis, such elements $x$ and $y$ exist in $A$. We then obtain a
contradiction as in \cite[p. 45]{C2}.   \end{proof}

\begin{proposition}
\label{algclsd2}
If $M$ is defined as in section 2, then $A$ is algebraically closed in $\operatorname{End}(M)$.
\end{proposition}

\begin{proof} For any $\sigma \in \operatorname{End}(M)$, there is a closed
unbounded subset $C_{\sigma }$ of $\omega _{1}$ such that for all $\alpha
\in C_{\sigma }$, $\sigma [M_{\alpha }]\subseteq M_{\alpha }$. For any $%
\sigma _{1},...,\sigma _{n}$ in $\operatorname{End}(M)$, choose $\alpha <\beta $
in $C_{\sigma _{1}}\cap ...\cap C_{\sigma _{n}}$ so that also $\alpha \in E$%
. Then each $\sigma _{i}$ induces an endomorphism, also denoted $\sigma _{i}$%
, of $M_{\beta }/M_{\alpha }.$ The endomorphism ring of $M_{\beta
}/M_{\alpha }$ is $\operatorname{End}$($G\oplus \mathbf{Z}^{(\omega )}$) and
restriction to $G$ defines a natural homomorphism, $\pi $, of $\operatorname{End}$%
($G\oplus \mathbf{Z}^{(\omega )}$) onto $\operatorname{End}(G)\cong A$ because $%
\operatorname{Hom}(G,\mathbf{Z}^{(\omega )})=0$. If $\sigma _{i}=a\in A$
(regarded as an element of $\operatorname{End}(M)$), then $\pi (a)=a$. Hence if $%
\sigma _{1},...,\sigma _{m}$ satisfy some ring equations over $A$, then so
do $\pi (\sigma _{1}),...,\pi (\sigma _{m})$.   \end{proof}

\smallskip\ 

Propositions \ref{algclsd} and \ref{algclsd2} provide an alternative proof
of ($\dag \dag $).

\

\end{document}